\documentclass[english,11pt]{article}
\usepackage[T1]{fontenc}
\usepackage{babel}
\usepackage{amsmath,amsfonts,amsthm}
\usepackage{color}
\usepackage{pdfsync}
\usepackage{graphicx}
\usepackage{url}

\DeclareMathSymbol{\leqslant}{\mathalpha}{AMSa}{"36} 
\DeclareMathSymbol{\geqslant}{\mathalpha}{AMSa}{"3E} 

\newcommand{\EE}{{\rm I\kern-2pt E}}
\newcommand{\RR}{{\rm I\kern-2pt R}}
\newcommand{\DD}{{\rm I\kern-2pt D}}
\newcommand{\PP}{{\rm I\kern-2pt P}}
\newcommand{\NN}{{\rm I\kern-2pt N}}

\begin{document}



\title{Am\'elioration de simulations Monte Carlo par des formes de Dirichlet}

\vspace{-2.6cm}


\title{Improving Monte Carlo simulations by Dirichlet forms}

\author{Nicolas Bouleau}
\date{---}
\maketitle



\noindent{\bf Abstract}

Equipping the probability space with a local Dirichlet form with square field operator $\Gamma$ and generator $A$
allows to improve Monte Carlo simulations of expectations and densities as soon as we are able to simulate a random variable $X$ 
together with $\Gamma[X]$ and $A[X]$. We give examples on the Wiener space, on the Poisson space and on the Monte Carlo space. When $X$ is real-valued
 we give an explicit formula yielding the density at the speed of the law of large numbers.



\noindent{\bf R\'esum\'e}

Nous montrons que, dans les situations o\`u l'espace de probabilit\'e est \'equip\'e d'une forme de Dirichlet locale
 avec carr\'e du champ $\Gamma$ et g\'en\'erateur $A$, la possibilit\'e de simuler une variable al\'eatoire $X$ ainsi que $\Gamma[X]$ et $A[X]$ 
permet d'acc\'el\'erer le calcul de l'esp\'erance de $X$ et de sa densit\'e. Nous donnons 
des exemples dans les cas de l'espace de Wiener,  de l'espace de Poisson et de l'espace de Monte Carlo. 
Lorsque $X$ est \`a valeurs r\'eelles
 nous donnons une formule explicite
 permettant d'obtenir la densit\'e \`a la vitesse de la loi des grands nombres. {\it Pour citer cet article~: N. Bouleau, C. R.
Acad. Sci. Paris, Ser. I ... (2005).}



\selectlanguage{english}


\section{Introduction}

  The efficiency of Dirichlet forms is known in order to obtain  existence of densities under weak hypotheses (cf [3]). 
We show here that they are still usefull for the computation of such densities.
Our framework is an error structure $(\Omega, {\mathcal A}, \PP,\DD,\Gamma)$, i.e. a probability space equipped 
with a local Dirichlet form $({\mathcal E},\DD)$ admitting a square field operator $\Gamma$ (cf [2],[3]). The associated 
$L^2$-generator is denoted $(A,{\mathcal D}A)$.

We consider a random variable $X\in{\mathcal D}A$ such that $X$, $\Gamma[X]$ and $A[X]$ are simulatable.\\

\noindent Example 1. Wiener space.

Let us consider a stochastic differential equation (sde) defined on the Wiener space equipped with the Ornstein-Uhlenbeck
error structure (cf [2],[3])
$$
X_t=x_0+\int_0^t \sigma(X_s,s)dB_s+\int_0^t r(X_s,s) ds
$$
By the functional calculus for the operators $\Gamma$ and $A$, if the coefficients are smooth, the triplet
 $(X_t, \Gamma[X_t], A[X_t])$ is a diffusion, solution to the equation
$$
\begin{array}{l}
\left(\begin{array}{l}
X_t\\
\Gamma[X_t]\\
A[X_t]
\end{array}\right)
=\left(\begin{array}{l}
x_0\\
0\\
0
\end{array}\right)
+{\displaystyle\int_0^t}
\left[
\begin{array}{ccc}
\sigma(X_s,s)&0&0\\
0&2\sigma^\prime_x(X_s,s)&0\\
-\frac{1}{2}\sigma(X_s,s)&\frac{1}{2}\sigma^{\prime\prime}_{x^2}(X_s,s)&\sigma'_x(X_s,s)
\end{array}
\right]
\left(\begin{array}{l}
1\\
\Gamma[X_s]\\
A[X_s]
\end{array}\right)
dB_s
\\
\\
\hspace{2.5cm}+{\displaystyle\int_0^t}
\left[
\begin{array}{ccc}
r(X_s,s)&0&0\\
\sigma^2(X_s,s)&2r^\prime_x(X_s,s)+\sigma^{\prime 2}_x(X_s,s)&0\\
0&\frac{1}{2}r^{\prime\prime}_{x^2}(X_s,s)&r^\prime_x(X_s,s)
\end{array}
\right]
\left(\begin{array}{l}
1\\
\Gamma[X_s]\\
A[X_s]
\end{array}\right)
ds
\end{array}
$$
Denoting $Y_t$ the column vector $(X_t,\Gamma[X_t],A[X_t])$ this equation writes
$
Y_t=Y_0+\int_0^ta(Y_s,s)dB_s+\int_0^tb(Y_s,s)ds
$
and applying the Euler scheme with mesh $\frac{1}{n}$ on [0,T] :
$
Y_t^n=Y_0+\int_0^ta(Y_{\frac{[ns]}{n}}^n,\frac{[ns]}{n})dB_s+\int_0^tb(Y_{\frac{[ns]}{n}}^n,\frac{[ns]}{n})ds.
$
 yields a process $Y_t^n=(X_t^n, (\Gamma[X])_t^n,(A[X])_t^n)^t$ for which it is easy to verify that
$
\Gamma[X_t^n]=(\Gamma[X])_t^n$ and
$A[X_t^n]=(A[X])_t^n$.

By known results (cf [1] [4] [5]) in order to compute the density of $X_T$,  we may approximate it 
by the solution $X_T^n$ of the Euler scheme. Thus, we have then to simulate $X_T^n$ in a situation where we are also able to simulate
$\Gamma[X_T^n]$ and $A[X_T^n]$.\\

\noindent Example 2. Poisson space.

Let $(\RR^d,{\mathcal B}(\RR^d), \mu,{\mbox{\bf d}},\gamma)$ be an error structure on $\RR^d$, $(a,{\mathcal D}a)$ its generator.
Let $N$ be a Poisson point process defined on $(\Omega, {\mathcal A}, \PP)$ with state space $\RR^d$ and
 intensity measure $\mu$. $(\Omega, {\mathcal A}, \PP)$ may be equipped with a 
so-called ``white" error structure $(\Omega, {\mathcal A}, \PP,\DD,\Gamma)$ (cf [2]) 
with the following properties : if $h\in{\mathcal D}a $ then $N(h)\in{\mathcal D}A $,  $\Gamma[N(h)]=N(\gamma[h])$
and $A[N(h)]=N(a[h]).$

In order to simulate $N(\xi)$ we have only to draw a finite (poissonian) number of i.i.d. random variables 
with law $\mu$ so that we are indeed in a situation where
$N(h)$, $\Gamma[N(h)]$, and $A[N(h)]$ are simulatable.\\

\noindent Example 3. Monte Carlo space.

Let $X=F(U_0,U_1,\ldots,U_m,\ldots;V_0,V_1,\ldots,V_n,\ldots)$ be a random variable defined on the space

\noindent$([0,1]^{\NN},{\mathcal B}([0,1]^{\NN}),dx^{\NN})\times([0,1]^{\NN},{\mathcal B}([0,1]^{\NN}),dx^{\NN})$ where
the $U_i$ are the coordinates of the first factor with respect to which $X$ is supposed to be regular, $V_j$ the ones of the second
factor with respect to which $X$ is supposed to be irregular or discontinuous (rejection method, etc.).
 
Let us put on the $U_i$ the following error structure
$$([0,1]^{\NN},{\mathcal B}([0,1]^{\NN}),dx^{\NN},\DD,\Gamma)= ([0,1],{\mathcal B}([0,1]),dx, {\mbox{\bf d}},\gamma)^{\NN}
$$
where $({\mbox{\bf d}},\gamma)$ is the closure of the operator $\gamma[u](x)=x^2(1-x)^2 u^{\prime2}(x)$ for $u\in{\mathcal C}^1([0,1]).$

Then under natural regularity assumptions, we have $\Gamma[X]= \sum_{i=0}^\infty F^{\prime2}_i U_i^2(1-U_i)^2.$ and
$$
A[X]=\sum_{i=0}^\infty(\frac{1}{2}F^{\prime\prime}_{ii} U_i^2(1-U_i)^2+F^\prime_i U_i(1-U_i)(1-2U_i))\\
$$
so that $X$, $\Gamma[X]$ and $A[X]$ are simulatable.
\section{Diminishing the bias}
Let  $(\Omega, {\mathcal A}, \PP,\DD,\Gamma)$ be an error structure.
For $X\in({\mathcal D}A)^d$,  $\underline{\underline{\mbox{var}}}[X]$ denotes the covariance matrix of $X$, $A[X]$ 
the column vector with components 
$(A[X_1],\ldots, A[X_d])$, $\underline{\underline{\Gamma}}[X]$ is the matrix $\Gamma[X_i,X_j]$ 
and $\sqrt{\underline{\underline{\Gamma}}[X]}$ denotes the positive symmetric square root of $\underline{\underline{\Gamma}}[X]$.

We follow the idea that the random variable $X+\varepsilon A[X]+\sqrt{\varepsilon}\sqrt{\underline{\underline{\Gamma}}[X]}\,G$
where $G$ is an exogeneous independent reduced Gaussian variable, has almost the same law as $X$. 
Starting from the fundamental relation of the functional calculus on $A$, 
an integration by parts argument gives the following lemma.\\

\noindent{\bf Lemma 2.1} {\it Let  $X\in({\mathcal D}A)^d$. we suppose that $X$ possesses a conditional density $\eta(x,\gamma,a)$
given $\underline{\underline{\Gamma}}[X]\!=\!\gamma$ et $A[X]\!=\!a$ such that $x\mapsto\eta(x,\gamma,a)$ 
be  ${\mathcal C}^2$ with bounded derivatives. Then $\forall x\in\RR^d$}
$$
\EE[-(A[X])^t\nabla_x\eta(x,\underline{\underline{\Gamma}}[X],A[X])
+\frac{1}{2}{\mbox{trace}}\left(\underline{\underline{\Gamma}}[X].\mbox{Hess}_x\eta\right)(x,\underline{\underline{\Gamma}}[X],A[X])]=0.$$

\noindent{\bf Theorem 2.2} {\it Let  $X$ be as in the preceding lemma, the conditional density $\eta(x,\gamma, a)$
 being ${\mathcal C}^3$ bounded with bounded derivatives. 
When $\varepsilon\rightarrow 0$, the quantity
$$\frac{1}{\varepsilon^2}\left(\EE[g(x-X-\varepsilon A[X], \varepsilon\underline{\underline{\Gamma}}[X])]-f(x)\right)$$
has a finite limit equal to }
$$
\frac{1}{2}\EE[(A[X])^t(\mbox{Hess}_x \eta)(x,\underline{\underline{\Gamma}}[X],A[X])A[X]
-\sum_{i,j,k=1}^dA[X_i]\Gamma[X_j,X_k]\eta^{\prime\prime\prime}_{x_ix_jx_k}(x,\underline{\underline{\Gamma}}[X],A[X])].
$$
{\it Proof.} If we write
$\;\;
\EE[g(x-X-\varepsilon A[X], \varepsilon\underline{\underline{\Gamma}}[X])]=
\int\mu(d\gamma,da)\int g(x-y-\varepsilon a, \varepsilon\gamma)
\eta(y,\gamma,a)dy 
$

\noindent
$
=\int\mu(d\gamma,da)\EE \eta(x-\varepsilon a-\sqrt{\varepsilon}\sqrt{\gamma}G,\gamma,a)
$
where $G$ is an $\RR^d$-valued reduced Gaussian variable, and if we expand with respect to $\sqrt{\varepsilon}$ and take the expectation,
terms in $\sqrt{\varepsilon}$ and $\varepsilon\sqrt{\varepsilon}$ vanish because $G$ and $G^3$ are centered and the term in 
$\varepsilon$ vanishes also thanks to the lemma. This gives the result.

About the variance, we obtain\\

\noindent{\bf Proposition 2.3}  {\it Let $X$ satisfying the assumptions of the lemma and such that} $(\mbox{det}
\underline{\underline{\Gamma}}[X])^{-\frac{1}{2}}\in L^1$, then
$$
\lim_{\varepsilon\rightarrow 0}\varepsilon^{d/2}\EE g^2(x-X-\varepsilon A[X],\varepsilon\underline{\underline{\Gamma}}[X])
=\lim_{\varepsilon\rightarrow 0} \varepsilon^{d/2}\mbox{var}g(x-X-\varepsilon A[X],\varepsilon\underline{\underline{\Gamma}}[X])
=\EE\left[\frac{\eta(x,\underline{\underline{\Gamma}}[X],A[X])}{(4\pi)^{d/2}\sqrt{{\mbox{\small det}}
\underline{\underline{\Gamma}}[X]}}\right].
$$

The quantity $\EE g(x-X-\varepsilon A[X],\varepsilon \underline{\underline{\Gamma}}[X])$ is obtained 
by simulation with the law of large numbers, so that the approximation $\hat{f}$ of the density $f$ of $X$ is
$$\hat{f}(x)=\frac{1}{N}\sum_{n=1}^N g(x-X_n-\varepsilon(A[X])_n,\varepsilon(\underline{\underline{\Gamma}}[X])_n)$$
where the indices $n$ denote independent drawings.
The preceding results show that, with respect to the usual kernel method, the speed, in the sense of the $L^2$-norm, is the same
as if the dimension was divided by 2.
\section{Direct formulae}
In the case where $X$ is real-valued, if in addition to $X$, $A[X]$, $\Gamma[X]$ we are able to simulate 
$\Gamma[X,\frac{1}{X}]$,  it is possible to obtain the dentity of $X$ at the speed of the law of large numbers thanks to the following
formulae :\\

\noindent{\bf Theorem 3.1} {\it a) If $X\in{\mathcal D}A$ with $\Gamma[X]\in\DD$ and $\Gamma[X]>0$ a.s. then $X$ 
has a density $f$
which possesses an  l.s.c. version $\tilde{f}$ given by }
$$
\tilde{f}(x)=\lim_{\varepsilon\downarrow 0}\uparrow\frac{1}{2}\EE\left(\mbox{sign}(x-X)
(\Gamma[X,\frac{1}{\varepsilon+\Gamma[X]}]+\frac{2A[X]}{\varepsilon+\Gamma[X]})\right).
$$
{\it b) 
If in addition $\frac{1}{\Gamma[x]}\in\DD$,
then $X$ has a density $f$ which is absolutely continuous and given by }
$$
f(x)=\frac{1}{2}\EE\left(\mbox{sign}(x-X)(\Gamma[X,\frac{1}{\Gamma[X]}]+\frac{2A[X]}{\Gamma[X]})\right).
$$

The proof is based on the relation
$$
\EE[\varphi^{\prime\prime}(X)\frac{\Gamma[X]}{\varepsilon+\Gamma[X]}]=
-\EE[\varphi^\prime(X)(\Gamma[X,\frac{1}{\varepsilon+\Gamma[X]}]+\frac{2A[X]}{\varepsilon+\Gamma[X]})].
$$
valid for any ${\mathcal C}^2$-function $\varphi$ with bounded derivatives which comes from the functional calculus 
using the general relation ${\mathcal E}[u,v]=-<A[u],v>$ $\forall u\in{\mathcal D}A$ $\forall v\in\DD$, and then 
applying it with $\varphi= \sqrt{\lambda^2+(y-x)^2}$ in order to get the monotone convergence result.\\

Under the hypotheses of theorem 3.1, as soon as $G\in\DD\cap L^\infty$, there are similar formulae for 
conditional expectations $\EE[G|X=x]$ : 
$$
f(x)\EE[G|X=x]=\frac{1}{2}\EE\left(\mbox{sign}(x-X)
(\Gamma[X,\frac{G}{\Gamma[X]}]+\frac{2GA[X]}{\Gamma[X]})\right).
$$
Let us finally remark that in these formulae, the factor on the right  of $\mbox{sign}(x-X)$ is centered and a variance optimisation may 
be performed thanks to an arbitrary deterministic function as done in [4] where direct formulae similar to those of section 3 are given in the case
of the Wiener space involving Skorokhod integrals instead of Dirichlet forms.

\end{document}